\documentclass{IEEEtran}
\usepackage{graphicx, algorithmic, algorithm}
\usepackage{paralist}
\usepackage{booktabs, multirow, fixmath, units}
\usepackage{psfrag, subfigure, bm}
\usepackage[usenames,dvipsnames]{pstricks}
\usepackage{amsmath,amsfonts,amssymb, mathabx}
\usepackage{cite}
\usepackage{pst-grad}
\usepackage[lowtilde]{url}
\def\x{\mathbold{x}}
\def\y{\mathbold{y}}
\def\v{\mathbold{v}}
\def\g{\mathbold{g}}

\def\W{\mathbf{W}}
\def\X{\mathbf{X}}

\def\myexp{\mathbb{E}}

\def\myset{}
\def\myom{\mathbold{\omega}}
\DeclareMathOperator*{\minimize}{\mathrm{minimize}}
\def\proj{{\mathsf{{P}}}}
\def\transp{\mathsf{T}}

\newtheorem{lemma}{Lemma}
\newtheorem{remark}{Remark}
\newtheorem{problem}{Problem}

\newtheorem{theorem}{Theorem}
\newtheorem{assumption}{Assumption}
\newenvironment{myproof}{\color{blue} \emph{Proof.} \ }{ \hfill $\blacksquare$ \color{black} }
\definecolor{mygray}{rgb}{0.9 0.85 0.85}

\renewcommand{\baselinestretch}{0.97}

\title{
Distributed Time-Varying Stochastic Optimization and Utility-based Communication} %\\ \texttt{Confidential Draft}

\author{Andrea Simonetto, Leon Kester, and Geert Leus  %
\thanks{A. Simonetto and G. Leus are with the Faculty of EEMCS, Delft University of Technology, 2826CD Delft, The Netherlands. e-mails: \footnotesize{$\{$a.simonetto, g.j.t.leus$\}$@tudelft.nl}. L. Kester is with \textsc{tno}, Oude Waalsdoorperweg 63, 2597AK The Hague, The Netherlands. email: \footnotesize{leon.kester@tno.nl}. This research was supported in part by STW under the D2S2 project from the ASSYS program (project 10561).}}%

\begin{document}

\maketitle

\begin{abstract}
We devise a distributed asynchronous stochastic $\epsilon$-gradient-based algorithm to enable a network of computing and communicating nodes to solve a constrained discrete-time time-varying stochastic convex optimization problem. Each node updates its own decision variable only once every discrete time step. Under some assumptions (among which, strong convexity, Lipschitz continuity of the gradient, persistent excitation), we prove the algorithm's asymptotic convergence in expectation to an error bound whose size is related to the constant stepsize choice $\alpha$, the variability in time of the optimization problem, and to the accuracy $\epsilon$. Moreover, the convergence rate is linear. Then, we show how to compute locally stochastic $\epsilon$-gradients that depend also on the time-varying noise probability density function (\textsc{pdf}) of the neighboring nodes, without requiring the neighbors to send such  \textsc{pdf}s at each time step. We devise utility-based policies to allow each node to decide whether to send or not the most up-to-date \textsc{pdf}, which guarantee a given user-specified error level $\epsilon$ in the computation of the stochastic $\epsilon$-gradient. Numerical simulations display the added value of the proposed approach and its relevance for estimation and control of time-varying processes and networked systems. %Our results are proven by using tools from linear algebra, notably bounds on the eigenvalues of Laplacian, convex analysis, and Peter-Paul inequalities. 
\end{abstract}

\section{Introduction}

We consider a time-varying stochastic optimization problem defined on time-varying functions that are distributed over a network of computing and communicating nodes. Let the nodes be labeled with $i\in\myset{V}=\{1,\dots,n\}$, and for each discrete time $k\in\mathbb{N}$, we equip each of them with the local function $f_{i,k}(\x, \myom): \mathbb{R}^{d\times q} \to \mathbb{R}$. In particular, the common vector $\x\in \mathbb{R}^d$ represents the decision variables, while the stochastic vector $\myom\in\Omega\subseteq\mathbb{R}^q$ is a stochastic variable drawn from a given (or estimated) time-varying probability density function (\textsc{pdf}) $p_{\myom,k}(\myom)$ on $\Omega$. We assume that $\myom$ is comprised of local stochastic variables as $\myom = (\myom_1^\transp, \dots, \myom_n^\transp)^\transp$, each of them of possible different dimensions and each of them uncorrelated with one another. Furthermore, we let $\myom_i \in \Omega_{i}$, $\Omega = \prod_{i=1}^n \Omega_{i}$, and $\myom_i \sim p_{\myom_i,k}(\myom_i)$ on $\Omega_{i}$. 

The main goal for the computing nodes at each discrete time $k$ is to solve the optimization problem 
\begin{equation}\label{optimalproblem}
\minimize_{\x \in\myset{X}_k} \myexp_{\myom,k} \hskip-0.05cm\Big[\hskip-0.05cm\sum_{i\in\myset{V}}\hskip-0.05cm f_{i,k}(\x, \myom)\Big]\hskip-0.1cm := \hskip-0.2cm \int_{\Omega} \sum_{i\in\myset{V}} \hskip-0.05cm f_{i,k}(\x, \myom) p_{\myom,k}(\myom) \mathrm{d}\myom,
\end{equation}
where each of the $f_{i,k}(\x, \myom)$ is a convex function of $\x$ for all $\myom \in \Omega$, while $\myset{X}_k$ is a compact convex set. 

We allow the computing nodes to communicate with their immediate neighbors defined via the undirected communication graph $\myset{G} = (\myset{V}, \myset{E})$, with edge set $\myset{E}$. In particular, each node $i$ can communicate with all the nodes $j\in \myset{N}_i:= \{j\in\myset{V}| (i,j)\in\myset{E} \}$. This communication possibility is not assumed to be synchronous among the computing nodes, but can happen at asynchronous times. We only assume that when node $i$ communicates with node $j$
, node $j$ communicates with node $i$ as well (that is, we assume an edge-asynchronous protocol).  

Stochastic optimization problems like~\eqref{optimalproblem} are rather popular in machine learning~\cite{Xiao2010, Duchi2013}. In this paper, we are more interested in their connection to the field of reasoning under uncertainty, self-awareness, and optimization-based design for distributed systems~\cite{vanFoeken2009, Ditzel2013}, where each computing node may learn or change the estimate of its local $p_{\myom_i,k}$ with time. In addition to that, other instances of the problem encompass distributed stochastic control \cite{Ueda2007}, finance \cite{Boyd2014}, and distributed statistical signal processing~\cite{Jakubiec2013}. As for possible high-impact societal applications, distributed time-varying stochastic problems as~\eqref{optimalproblem} arise in the context of smart grids~\cite{Zhao2012, Low2014} and cooperative adaptive cruise control~\cite{PATH2008}, to name a few.

If, for the moment, we focus only at the time-varying nature of the problem, we can also recognize in~\eqref{optimalproblem} optimization programs that appear in distributed estimation of stochastic time-varying signals~\cite{Jakubiec2013}, in distributed control of mobile multi-robot systems with time-varying tasks~\cite{Tu2011}, and as a result of sequential convex programming approaches to multi-agent non-convex problems~\cite{Simonetto2012a}. When each of the functions $f_{i,k}(\x)$ and the set $\myset{X}_k$ are time-invariant, several approaches can be applied to solve~\eqref{optimalproblem}. These techniques differ for the assumptions they require and the properties they can ensure (convergence, convergence rate, resilience to asynchronous communication protocols, among others). Examples of such approaches are the stochastic subgradient~\cite{Srivastava2011}, dual averaging~\cite{Duchi2012}, and the alternating direction method of multipliers~\cite{Boyd2011}. %The latter technique has been shown to be the most general with respect to the assumptions it imposes on the original problem~\eqref{optimalproblem} to guarantee convergence. Furthermore, it results in a linear convergence rate, which is considered to be fast in comparison with other existing techniques. 
Since the aforementioned techniques are iterative, and they require communication among the nodes to converge to an optimizer of~\eqref{optimalproblem}, they would \emph{provably} converge in the case of time-varying $f_{i,k}(\x)$'s and/or time-varying $\myset{X}_k$ only when each node could exchange an infinite number of messages with its neighbors, between  consecutive time steps $k$ and $k+1$. Specific methods that account for a finite number of messages between consecutive time steps and still guarantee convergence have been proposed in~\cite{Kamgarpour2008, Braca2010,Cattivelli2010a, Farina2010, Tu2011, Bajovic2011, Zavlanos2013, Jakubiec2013, Cavalcante2013, Ling2013}, but they are all limited to specific \emph{deterministic} versions of~\eqref{optimalproblem}. Notably, in~\cite{Jakubiec2013, Ling2013}, the authors work under the same assumptions that we will use, however they consider deterministic and \emph{unconstrained} optimization problems, while in~\cite{Cavalcante2013}, the authors employ deterministic subgradient methods and assume that the optimizers of~\eqref{optimalproblem} do not change in time. 

\textbf{Contributions.} As a first contribution, we consider the stochastic nature of the optimization problem and we propose an asynchronous stochastic gradient-based distributed algorithm for the computing nodes to converge to an optimizer of~\eqref{optimalproblem}. In fact, due to the time-varying nature of the problem and due to possible errors in the computation of the stochastic gradient, the convergence will be shown up to an error bound, whose size is directly dependent on the mentioned elements. This algorithm can be seen as a generalization in a time-varying context of the work in~\cite{Srivastava2011} where only one iteration of the algorithm is performed between consecutive time steps, as well as a generalization of the work in~\cite{Jakubiec2013, Ling2013} in constrained, asynchronous, and stochastic settings. In addition, in contrast to~\cite{Jakubiec2013, Ling2013}, our algorithm does not hinge on dual variables to reach a common decision vector among the nodes (which complicates significantly the theoretical analysis of convergence), but is instead based on consensus protocols, which are easier to analyze and embed on real hardware.   

The proposed algorithm can also be seen as a stand-alone contribution, that can be applied to many different scenarios, for example, in deterministic settings or when the nodes can compute the gradient exactly (we will show an example of these cases in the numerical simulations).  

As a second contribution, we consider the case in which the computing nodes have access to the local \textsc{pdf}s $p_{\myom_i,k}$ only. In this case, to avoid the communication overhead of sending at  each time step $k$, the most current local \textsc{pdf} to the other nodes, we devise a utility-based policy that enables each node to decide locally whether to send such $p_{\myom_i,k}$ or not. This utility-based policy guarantees that despite the possible outdated information on $p_{\myom_i,k}$, each of the nodes is able to compute their stochastic gradient up to a user-specified precision $\epsilon$. This fits perfectly in the scheme of the proposed asynchronous distributed algorithm, and is  particularly useful to limit the overhead of sending at each time step the varying \textsc{pdf}s, which might be non-Gaussian and difficult to encode in a message with a limited number of parameters. This second contribution can be seen as a generalization of event-triggered optimization~\cite{Wan2009}, where we determine the triggering mechanism not only based on convergence arguments, but also 
on performance guarantees.

Interesting related work for event-based mechanisms can be found in event-triggered control~\cite{Zhong2010, Heemels2012}, estimation~\cite{Trimpe2012}, and energy constrained communication networks~\cite{Arroyo-Valles2009, Razavi2013}. An important difference between our contribution and these fields of research is that we cast the global optimization problem as a stochastic program, which is able to capture \emph{unpredicted (i.e., non-modeled) but locally measurable changes} and judge \emph{locally} their importance to the \emph{global} performance.  

%Finally, as a third contribution, we propose explicit trade-offs between communication and computation costs for the nodes. These trade-offs are based on upper bounds on the utility-based policy, which allow the nodes to decide whether to send the most up-to-dated version of $p_{\myom_i,k}$ with reduced computation cost, yet on average more frequently. In the limit of sending with no-computation, we will prove that the utility-based policy leads to a periodic communication mechanism, as happens often in event-based protocols~\cite{Trimpe2012}. 

\textbf{Organization.} The first contribution of the paper is encoded in Problem~\ref{p1}, which is formally expressed in Section~\ref{sec:p1} and whose solution is outlined in Algorithm~\ref{alg} and characterized in Theorem~\ref{convergence}. The second contribution is the approach to tackle Problem~\ref{p2}, which is fully discussed in Section~\ref{sec:p2} and Theorem~\ref{policy}. %Section~\ref{sec:c3} is instead devoted to the third contribution of the paper. 
Numerical results support the theoretical findings and are reported in Section~\ref{sec:num}. All the proofs are grouped in Section~\ref{proofs}.

 %The remainder of the paper is organized as follows. In Section~\ref{sec:as}, we present the main assumptions and a formal problem formulation. Sections~\ref{sec:p1} and \ref{sec:p2} are devoted to the discussion and analysis of the proposed approach, while in Section~\ref{sec:num}, a numerical example is shown to assess the solution in practice. Finally, in Section~\ref{sec:con}, we draw our conclusions and we define the foreseeable next steps. Proofs of the theorems are given in Appendix.   

\textbf{Notation.} For any vector $\x\in\mathbb{R}^n$, the norm $\|\x\|$ represents the standard Euclidean norm. For any real-valued squared matrix $\X\in\mathbb{R}^{n\times n}$, we say $\X \succeq 0 $ or $\X \preceq 0$ if the matrix is positive semi-definite or negative semi-definite, respectively. For any positive semi-definite matrix of dimension $n$, we order its eigenvalues as $0\leq \lambda_1(\X) \leq \dots \leq \lambda_n(\X)$. The symbol $\otimes$ represents the Kronecker product. Given a
 differentiable convex function $f(\x): \myset{X}\subseteq \mathbb{R}^n \to \mathbb{R}$ and a non-negative scalar $\epsilon$, an $\epsilon$-(sub)gradient of $f(\x)$ at $\x\in X$ is a vector $\tilde{\g}\in\mathbb{R}^n$ such that
\begin{equation}\label{egrad}
\tilde{\g}^\transp (\y-\x) \leq f(\y) - f(\x) + \epsilon, \quad \textrm{for all } {\y} \in \myset{X}
,
\end{equation}
and finally, given a compact set $X$, the symbol $|X|$ stays for $|X|:= \max_{\x\in X}\{||\x||\}$.  
%holds almost surely. The chosen $\tilde{\g}_{i,k}$ has the form

\section{Assumptions and Problem Statement}\label{sec:as}

Let $X_{k,k-1}$ be defined as the convex hull of the union of the sets $X_k$ and $X_{k-1}$, for $k\geq 1$, i.e., $X_{k,k-1}:= \mathrm{cvxh}\{X_k \cup X_{k-1}\}$, which will be used in our proofs. Throughout the text we assume the following simplifying assumptions. 

\begin{assumption}\label{as.function}\emph{(Sets and Objective functions)}
The sets $X_k$ are compact convex sets for all $k\geq 0$. The functions $f_{i,k}(\x, \myom)$, $k\geq 1$, are twice differentiable and convex for all $\x \in X_{k,k-1}$, $\myom \in \Omega$, and they have the following properties: 
\vskip0.1cm
\hspace*{-0.5cm}\begin{minipage}{0.495\textwidth}
\begin{enumerate}%[leftmargin=.55cm]
\item the expectations $\myexp_{\myom,k}[f_{i,k}(\x, \myom)]$ are strongly convex with respect to $\x \in X_{k,k-1}$, $k\geq 1$. In particular, the eigenvalues of the Hessian $\myexp_{\myom,k}[\nabla^2_{\x\x} f_{i,k}(\x, \myom)]$ are lower bounded by the strong convexity constant $m_f$ for all $\x \in X_{k,k-1}$,
$$
\myexp_{\myom,k}\big[\nabla^2_{\x\x} f_{i,k}(\x, \myom)\big] \succeq m_f {\bf I},
$$  
or equivalently
\begin{multline*}
\myexp_{\myom,k}\big[f_{i,k}(\y, \myom)-f_{i,k}(\x, \myom)\big] \geq \\\myexp_{\myom,k}\big[\nabla_\x  f_{i,k}(\x, \myom)^\transp\big] (\y-\x) + \frac{m_f}{2}\|\y-\x\|^2,
\end{multline*}
for all vectors $\x, \y$ in $X_{k,k-1}$;
\item the gradients of the expectations $\myexp_{\myom,k}[f_{i,k}(\x, \myom)]$ with respect to $\x$ are Lipschitz continuous with constant $L$, meaning that 
\begin{equation*}
\|\myexp_{\myom,k}\big[\nabla_{\x} f_{i,k}(\x, \myom) - \nabla_{\y} f_{i,k}(\y, \myom)\big]\|\leq L \|\x-\y\|, 
\end{equation*}
for all vectors $\x, \y$ in $X_{k,k-1}$, or equivalently
$$
\myexp_{\myom,k}\big[\nabla^2_{\x\x} f_{i,k}(\x, \myom)\big] \preceq L {\bf I},
$$
for all $\x \in X_{k,k-1}$.
\end{enumerate}\end{minipage}
\end{assumption}\vskip0.1cm

\begin{assumption}\label{as.stoc}\emph{(Stochastic variables)}
The stochastic variables $\myom_i$ affect the objective functions only locally, i.e., 
\begin{equation*}
\nabla_{\myom_j} f_{i,k}(\x, \myom) = {\bf 0},\quad \mathrm{ if~}  j\notin \myset{N}_i\cup i;
\end{equation*}
the local cost functions are separable w.r.t. the $\myom_i$, i.e., 
\begin{equation*}
f_{i,k}(\x, \myom) = \sum_{j\in\myset{N}_i\cup i} f_{i,k,j}(\x, \myom_j), \qquad \mathrm{for~} i=1,\dots,n;
\end{equation*}
and finally, each of the \textsc{pdf}s $p_{\myom,k}$ are defined over a compact set, that is to say that each of the $\Omega_{i}$ is a compact set. 
\end{assumption}

\begin{assumption}\label{as.opt}\emph{(Optimal problem)}
There exists a solution $\x^*_{k}$ for~\eqref{optimalproblem} at each time step $k$, and the distance between these solutions at two subsequent time steps is upper bounded as
\begin{equation*}
\|\x^*_{k+1} - \x^*_k\| \leq \delta_\x, \quad \textrm{for } k\geq 1,
\end{equation*}
where $\delta_\x$ is a finite non-negative scalar. 
\end{assumption}

\begin{assumption}\label{as.persistent}\emph{(Communication graph)} At each iteration $k$ the symmetric adjacency matrix of the communication graph $\mathbf{A}_{k}$ is generated by an i.i.d. Bernoulli process with $\mathrm{Pr}[[\mathbf{A}_k]_{ij} = 1]=s_{ij}>0$ for all $(i,j)\in\myset{E}$, with a given probability $0< s_{ij}\leq 1$. Furthermore, let $\myset{E}_k$ be the communication edge set at time $k$, and let $\myset{G}_k :=(\myset{V},\myset{E}_k)$. For every $k'\geq 1$, there exists an integer $T\geq0$ such that: 
\begin{enumerate}
\item the union of the edge sets satisfies $\bigcup_{\ell=k'}^{k'+T} \myset{E}_\ell = \myset{E} $;
\item the union graph, i.e., $\bigcup_{\ell=k'}^{k'+T}  \myset{G}_\ell$, is connected. 
\end{enumerate}
\end{assumption}

Assumptions~\ref{as.function} and \ref{as.opt} are typical in time-varying optimization problems; in particular, Assumption~\ref{as.function} ensures that the optimizer of \eqref{optimalproblem} (if it exists) is unique, while Assumption~\ref{as.opt} ensures the existence of said optimizers and bounds their variability, which will be crucial in the convergence of the proposed solution. For a broader perspective on these two assumptions, the reader is referred to~\cite[Chapter 6]{Polyak1987}. We notice the need of considering a larger set ($\myset{X}_k\cup\myset{X}_{k-1}$) for the validity of the smoothness conditions on the cost functions; this is due to the time-varying nature of the problem and its role will be clear in the proofs. The compactness assumption on $X_k$ is also quite used in ($\epsilon$)-subgradient algorithms and it is not very restrictive (one could always use a box constraint to limit the variability of the decision variable $\x$); furthermore, this assumption is not necessary when the computation of the gradient is exact. %, which is perhaps one of the first comprhensive treatment of time-varying optimization problems. 
Assumption~\ref{as.persistent} is typically required in asynchronous distributed protocols and is rather weak. Assumption~\ref{as.stoc} is a simplifying assumption on the nature of the stochastic interdependences; in particular Assumption~\ref{as.stoc} says that only the neighboring stochastic variables have an effect on each node's cost function $f_{i,k}$, and that this effect is the sum of local components coming from different neighbors. This is certainly quite reasonable in distributed systems. Finally, Assumption~\ref{as.stoc} says that the probability space is compact, which is not a ver
y restrictive requirement in practice and it is important in quantifying how important are the variations in $f_{i,k}$ w.r.t. the variations in $p_{\myom,k}$. 

As mentioned in the Introduction, the problem we are interested in is twofold. First of all, we want to enable the computing nodes to solve~\eqref{optimalproblem} in a distributed fashion, where each of the nodes communicates with their neighbors only. For this first task, we introduce local copies of the decision variable $\x_k$. These local copies are referred to as $\y_{i,k}$. We formally formulate the problem at hand as 

\begin{problem} \label{p1}
Devising a stochastic $\epsilon$-(sub)gradient distributed algorithm in order to enforce that the local decision variable $\y_{i,k}$ eventually converges (up to a bounded error) to the optimal solution of~\eqref{optimalproblem} at time step $k$ ($\x_{k}^*$), or formally,  
$$
\liminf_{k\to \infty} \myexp\left[\|\y_{i,k} - \x^*_k  \|^2\right] \leq \delta, \textrm{ for } i\in\myset{V},
$$
for some $\delta\geq0$, which has to decrease if $\delta_\x \to 0$ and $\epsilon \to 0$. 
\end{problem}

%This also implies convergence of the $\y_{i,k}$ to a common value up to a bounded error. This bounded error will comprise two terms; one that comes from possible communication/computation errors (the $\epsilon$ part in the gradient computation), say $\delta_\epsilon$, and one that comes from the time-varying nature of the problem, say $\delta_{\textrm{TV}}$. 

The second aspect of the problem we consider is to limit the need of communicating the most actual \textsc{pdf} $p_{\myom,k}$ among the nodes and can be formulated as 

\begin{problem} \label{p2}
Devising a utility-based policy that allows the nodes to decide whether to send or not their $p_{\myom_i,k}$ and yet that guarantees a prescribed level of accuracy $\epsilon$ in the computation of the stochastic $\epsilon$-(sub)gradient. 
\end{problem}

\section{Proposed Approach for Problem~\ref{p1}}\label{sec:p1}

In order to solve the first part of the problem we propose a distributed asynchronous stochastic gradient algorithm. We start by defining the time-varying matrix $\W_k$ and two different stepsizes $\alpha>0$ and $\beta>0$. The matrix $\W_k$ is a symmetric matrix based on the adjacency matrix ${\bf A}_k$ as,
\begin{equation}
[\W_k]_{i,j} = \left\{\begin{array}{cc} - [{\bf A}_k]_{i,j} & \textrm{for } j \neq i \\
\sum_{l=1}^n [{\bf A}_k]_{i,l}& \textrm{for } j = i \end{array}\right. .
\end{equation}
From Assumption~\ref{as.persistent}, the symmetric matrix $\W_k$ has nonzero elements if and only if
 the related nodes can communicate with each other, it is rank deficient and in particular $\W_k {\b
f 1}_n = {\bf 0}_n$, and finally, for the sequence of matrices $\{\W_k\}$, 
\begin{equation}
\myexp[\W_k] = \bar{\W} = \bar{\W}^\transp, \textrm{ with } \lambda_2(\bar{\W})> 0.  
\end{equation}

As for the stepsizes $\alpha$ and $\beta$, the former caters for convergence to the optimal solution of~\eqref{optimalproblem}, while the latter dictates the consensus among the different computing nodes. The proposed algorithm is outlined in Algorithm~\ref{alg}. 

\begin{algorithm}
\caption{Asynchronous distributed stochastic $\epsilon$-gradient}\label{alg}
\footnotesize
\begin{subequations}
\hskip .5cm Initialize by picking locally an arbitrary $\y_{i,1} \in \myset{X}_0$. Then for $k\geq 1$:

\begin{enumerate}
\item compute the local variable $\v_{i,k+1}$ by local communication as 
\begin{equation}\label{v}
\v_{i,k+1} = \y_{i,k} - \beta \sum_{j = 1}^n [\W_k]_{i,j} \y_{j,k}\,,
\end{equation}
for which we will show that in fact $\v_{i,k+1}\in\myset{X}_{k-1}$ in Lemma~\ref{lemma1};

\item compute locally the stochastic $\epsilon$-(sub)gradient of $f_{i,k}$ with respect to $\x$ at $\v_{i,k+1}$, as
\begin{eqnarray}\label{gtilde}
\tilde{\g}_{i,k} &=& \myexp_{{\myom'}}[\nabla_\x f_{i,k}(\x,{\myom})|_{\v_{i,k+1}}] := \\
&& \int_{\Omega} \left(\nabla_\x f_{i,k}(\x,{\myom})|_{\v_{i,k+1}}\right) p_{\myom'}(\myom) \mathrm{d}\myom, \nonumber 
\end{eqnarray}
for which $p_{\myom'}\in\Omega$ is the \textsc{pdf} of the stochastic variable $\myom$ know by node $i$ at time $k$ (which is possibly outdated). We give a sufficient condition for the chosen $\tilde{\g}_{i,k}$ to be a stochastic $\epsilon$-(sub)gradient in Lemma~\ref{lemma1};   

\item update the local variable $\y_{i,k}$ as
\begin{equation}\label{update}
\y_{i,k+1} = \proj_{\myset{X}_{k}}\left[\v_{i,k+1} - \alpha \tilde{\g}_{i,k}\right];
\end{equation}
where $\proj_{\myset{X}_k}[\cdot]$ indicates the projection operator; 

\item go to step 1.

\end{enumerate}
\end{subequations}
\end{algorithm}

The proposed algorithm involves only local communication (in step 1) and local computations. This is true, provided that each node can obtain the \textsc{pdf} $p_\myom'$ necessary in step 2. We show later in this section how this can be achieved by a utility-based local communication mechanism, that can be completely disconnected by the communication of the local $\y_{i,k}$'s. 

Our first result establishes convergence of the sequence of local variables $\{\y_{i,k}\}$ generated by Algorithm~\ref{alg} to the optimizer of~\eqref{optimalproblem} up to a bounded error. The error is comprised of two parts, one coming from the error in the computation of the gradient, the other coming from the time-varying nature of the original problem. 

As for Assumptions~\ref{as.function} and \ref{as.opt}, the gradient of the expectations at optimality are bounded for all $k$, and in particular,  
\begin{equation}\label{b.G}
\max_{i\in\myset{V}}\big\{\big\|\sum_{j\in\myset{V}, j\neq i} \myexp_{\myom,k}[\nabla_\x f_{j,k}(\x, \myom)|_{\x^*_k}]\big\|\big\} \leq G,
\end{equation}
for a certain nonnegative bound $G$.

\noindent\colorbox{mygray}{\begin{minipage}{0.48\textwidth}\begin{theorem}\label{convergence}
Let $\{\y_{i,k}\}$ be the sequence generated by Algorithm~\ref{alg} where each $\tilde{\g}_{i,k}$ is a stochastic $\epsilon$-(sub)gradient of $f_{i,k}$ w.r.t. $\x$ at $\v_{i,k+1}$ computed as~\eqref{gtilde}; let ${\g}_{i,k}$ be the up-to-date stochastic gradient as defined in~\eqref{gom}; let $\y_k$ be the stacked $nd$-dimensional column vector of the local variables $\y_{i,k}$. Let $G$ be defined as~\eqref{b.G}. %Let $\underline{s}$ be the minimum among the edge probabilities $0<s_{ij}\leq 1$, as defined in Assumption~\ref{as.persistent}. 
Define $ \varrho := 1 + \alpha^2 L^2- \alpha m_f $. Under Assumptions~\ref{as.function} till \ref{as.persistent}, if $\|\tilde{\g}_{i,k}-{\g}_{i,k}\|\leq \epsilon/2|X_{k,k-1}|$ for each $k\geq 1$, then by choosing $\beta~<~1/n$, and $\alpha~<~m_f/L^2$, the sequence $\{\y_{k}\}$ convergences as 
\begin{multline*}
\liminf_{k\to \infty}\myexp[\|\y_{k+1} - {\bf{1}}_n\otimes\x_{k+1}^{*}\|^2]\leq \\
\frac{1}{1-\varrho} \left(\alpha\hat{\psi}(\alpha, \epsilon, \sqrt{\gamma})\frac{1}{\sqrt{\gamma}} + n\delta_\x^2 \frac{1}{1-\sqrt{\gamma}}\right),
\end{multline*}
where $\gamma = 1 - \beta\lambda_2(\bar{\W})$ and 
$$
\hat{\psi}(\alpha, \epsilon,\sqrt{\gamma}) = n \frac{\epsilon}{\sqrt{\gamma}} \left(\alpha\frac{\epsilon}{4 |\myset{X}_{k,k-1}|^2} + 2 \alpha L  + 2 \right) + \frac{\alpha\, n\, G^2}{1-\sqrt{\gamma}}. 
$$ 
Furthermore, $0< \varrho <1$ and the convergence rate is linear. 
\end{theorem}\end{minipage}}

\vskip0.1cm As we will see in Section~\ref{proofs}, the proof of Theorem~\ref{convergence} is based on the Peter-Paul inequality for a particular choice of the parameter $\mu$. This choice is valid if $0<\gamma<1$, as in our case, and allows us to separate the effect of $\alpha$ with the effect of the network connectivity $\lambda_2(\bar{\W})$. Furthermore, it gives us some intuition on which errors will play a more predominant role for different $\lambda_2(\bar{\W})$. In particular, for highly connected graphs, when $\gamma \to 0$, then the term $\hat{\psi}(\alpha, \epsilon, \sqrt{\gamma})$ is more important than the time-varying term. This is expected since consensus mixing is faster and the nodes quickly agree on a common approximate optimizer. Furthermore, in this case, to reduce the error, we will have to choose smaller and smaller $\alpha$'s. On the contrary for poorly connected graphs, $\gamma \to 1$, the error term with $G$ and the time-varying one are dominant. A similar results was found in~\cite{Jakubiec2013} for deterministic strongly convex optimization problems using a \emph{synchronized} dual decomposition approach. 

%A similar results \emph{in terms of decomposition of the competences} for $\alpha$ and $\beta$ was found in~\cite{Simonetto2012b}. In particular, $\alpha$ is determined based on the function properties, while $\beta$ is based on the graph topology. 

\begin{remark}\label{tight}
The bound that has been derived is not expected to be tight for the whole range of $\gamma$'s; however one can use the same Peter-Paul inequality with a different choice of $\mu$ to obtain tighter bounds for specific values of $\gamma$. Our choice was determined based on the idea of dividing the terms of convergence ($\alpha$) and consensus ($\beta$). In addition, it is expected that the smaller $\delta_\x$ and $\epsilon$ are, the tighter the bound is.   
%Both $\alpha$ and $\beta$ might be computed locally if needed, in particular, by assuming that each node knows a local $\alpha_i$ linked with its own $f_{i,k}$, and each communication link is equiprobable, i.e., $s_{ij} = s$. Then, $\alpha$ can be computed by running a $\min$-consensus on the graph, while $\beta$, by lower bounding $\underbar{s}/\lambda_n(\bar{\W})$ as $1/n$ \cite{Rojo2007}.  
\end{remark}

\section{Proposed Approach for Problem~\ref{p2}}\label{sec:p2}

We focus now on the second part of our problem, which is how to generate $\epsilon$-(sub)gradient vectors for a specific value of $\epsilon$, that is, under the sufficient condition of Lemma~\ref{lemma1} in Section~\ref{proofs}, how to generate $\tilde{\g}_{i,k}$ such that $\|\tilde{\g}_{i,k} - \g_{i,k}\| \leq \epsilon/2 |X_{k,k-1}|$. In particular, we consider a utility-based event-triggered mechanism that allows the computing nodes to send one another the latest updated \textsc{pdf} of $\myom$ only when needed, to guarantee a prescribed level of $\epsilon$.

The error in computing the stochastic gradient at node $i$ is
\begin{multline*}   
{\g}_{i,k} - \tilde{\g}_{i,k} = \\ \myexp_{\myom,k} [\nabla_{\x} f_{i,k}(\x,\myom)|_{\v_{i,k+1} }] - \myexp_{\myom'} [\nabla_{\x} f_{i,k}(\x,\myom)|_{\v_{i,k+1}}] = \\ \hskip-0.2cm\sum_{\hskip.1cm j\in\myset{N}_i}\hskip-0.1cm\myexp_{\myom_j\hskip-0.05cm,k} [\nabla_{\x} f_{i,k,j}(\x,\myom_j)|_{\v_{i,k+1}}\hskip-0.05cm] - \myexp_{\myom_j'}\hskip-0.05cm[\nabla_{\x} f_{i,k,j}(\x,\myom_j)|_{\v_{i,k+1}}\hskip-0.05cm]\end{multline*}
and as a consequence of Lemma~\ref{lemma1}, we can derive the sufficient condition for each $j\in\myset{N}_i$
%\begin{multline*}   
%\left\|\myexp_{\myom,k} [\nabla_{\x} f_{i,k}(\x,\myom)|_{\v_{i,k+1} }] - \myexp_{\myom'} [\nabla_{\x} f_{i,k}(\x,\myom)|_{\v_{i,k+1}}] \right\| = \\ \left\|\sum_{\hskip1cm j\in\myset{N}_i\cup i}\left(\myexp_{\myom_j,k} [\nabla_{\x} f_{i,k,j}(\x,\myom_j)|_{\v_{i,k+1} }] - \myexp_{\myom_j'} [\nabla_{\x} f_{i,k,j}(\x,\myom_j)|_{\v_{i,k+1}}]\right) \right\|  \leq \frac{\epsilon_{i,k}}{2 |X_k|},
%\end{multline*}
%which indicates the sufficient condition for each $j\in\myset{N}_i$ 
\begin{multline}\label{cond}   
\hskip-.4cm \left\|\myexp_{\myom_j\hskip-0.05cm,k} [\nabla_{\x} f_{i,k,j}(\x,\myom_j)|_{\v_{i,k+1}}\hskip-0.05cm] - \myexp_{\myom_j'}\hskip-0.05cm[\nabla_{\x} f_{i,k,j}(\x,\myom_j)|_{\v_{i,k+1}}\hskip-0.05cm]\right\| \\ \leq  \frac{\epsilon}{2 |X_{k,k-1}| \Delta_i},
\end{multline}
where $\Delta_i$ is the number of possible neighbors of $i$, as in $\myset{E}$. 

\begin{figure}
\centering
\psfrag{i}{$i$}\psfrag{j}{$j$}
\psfrag{c}{$\pi_{ij} \to$}
\psfrag{d}{$\leftarrow $ $\pi_{ji}$}
\psfrag{a}{\footnotesize $\nabla_\x f_{i,k}(\x,\myom), \, \forall k $}
\psfrag{b}{\footnotesize $\nabla_\x f_{j,i}'(\x,\myom_i)|_{\v_j'}$}
\psfrag{e}{\footnotesize $\v_{i,k+1},\, \forall k$}
\psfrag{f}{\footnotesize $p_{\myom_i,k},\, \forall k$}
\psfrag{g}{\footnotesize $p_{\myom_j'}, \v_j'$}
\psfrag{h}{\footnotesize $\nabla_\x f_{j,k}(\x,\myom), \, \forall k$}
\psfrag{l}{\footnotesize $\nabla_\x f_{i,j}'(\x,\myom_j)|_{\v_i'}$}
\psfrag{m}{\footnotesize $\v_{j,k+1}, \, \forall k$}
\psfrag{n}{\footnotesize $p_{\myom_j,k}, \, \forall k$}
\psfrag{o}{\footnotesize $p_{\myom_i'}, \v_{i}'$}
\includegraphics[width=0.45\textwidth]{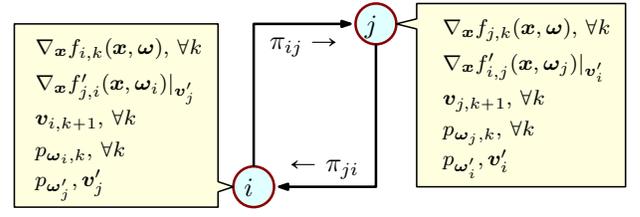}
\caption{Quantities known at each node and the communication policy $\pi_{ij}$, the quantities with 
$'$ are the outdated ones.}
\label{ut}
\end{figure}

If we could use the condition~\eqref{cond} to establish a sending policy, then $\tilde{\g}_{i,k}$ would be guaranteed to be an $\epsilon$-(sub)gradient and we would know the error floor in Theorem~\ref{convergence}. Unfortunately, neither node $i$ nor node $j$ can compute, check, or enforce~\eqref{cond}. This is because nodes have either their updated $\v_{i,k+1}$ or the updated $p_{\myom_j,k}$, b
ut not both at the same time (see Figure~\ref{ut} for a pictorial representation). However, along the mechanism to send the most up-to-date $p_{\myom_j,k}$, we can think of a similar mechanism to send or not the most up-to-date gradient $\nabla_{\x} f_{i,k,j}(\x,\myom_j)|_{\v_{i,k+1}}$, which is a function of $\myom_j$. Let $\nabla_{\x} f_{i,j}'(\x,\myom_j)|_{\v_{i}'}$ be the possibly outdated gradient function of node $i$ available at node $j$. We define the utility metrics, 
%\begin{align}
\begin{subequations}\label{metrics}
\begin{multline}
\mathcal{U}_{\mathrm{S},ji}^{(1)}(\myom_k|\myom') : = \\ \left\|\myexp_{\myom_j,k} [\nabla_{\x} f_{i,j}'(\x,\myom_j)|_{\v_{i}' }] - \myexp_{\myom_j'} [\nabla_{\x} f_{i,j}'(\x,\myom_j)|_{\v_{i}'}]\right\|,
\end{multline}
%\end{equation}
\begin{equation}
\hskip-2cm \mathcal{U}_{\mathrm{S},j}^{(2)}(\myom_k|\myom') : = \max_{\myom_j \in\Omega_j} |p_{\myom_j,k}(\myom) - p_{\myom'_j}(\myom)|,
\end{equation}
%\begin{equation}
\begin{multline}
\mathcal{U}_{\mathrm{R},ij}(\nabla_{\x} f_{i,k,j}|\nabla_{\x} f'_{i,j}) : = \\ \int_{\Omega_j}\Big\| \nabla_{\x} f_{i,k,j}(\x,\myom_j)|_{\v_{i,k+1}}- \nabla_{\x} f_{i,j}'(\x,\myom_j)|_{\v_{i}'} \Big\| \mathrm{d} \myom_j .
\end{multline}
%\end{align}
\end{subequations}

These utility functions specify how useful is the up-to-date information with respect to the outdated one. All of them are locally computable by nodes $i$ and $j$, respectively (Figure~\ref{ut}). With these quantities, we specify the sending policy $\pi_{ij}$ of node $i$ for each neighboring node $j$, for any scalar $\eta\in[0,1]$ and $\nu > 0$ as follows:
\begin{equation}\label{pi}
\pi_{ij} : \left\{\begin{array}{ll} \textrm{iff } \mathcal{U}_{\mathrm{R},ij}(\nabla_{\x} f_{i,k,j}|\nabla_{\x} f'_{i,j}) > \displaystyle\frac{\eta}{\nu}\frac{\epsilon}{2 |X_{k,k-1}| \Delta_i} ,& \\  \hskip3cm \textrm{send }  \nabla_{\x} f_{i,k,j}(\x,\myom_j)|_{\v_{i,k+1}} \textrm{ to } j & \\ & \\
 \textrm{iff } \left[\begin{array}{ll} \textrm{either} & \mathcal{U}_{\mathrm{S},ij}^{(1)}(\myom_k|\myom') > (1-\eta)\displaystyle\frac{\epsilon}{2 |X_{k,k-1}| \Delta_j}\\ \textrm{or} & \mathcal{U}_{\mathrm{S},i}^{(2)}(\myom_k|\myom') > \nu \end{array}\right] ,& \\ \hskip3cm \textrm{send the \textsc{pdf} } p_{\myom_i,k}   \textrm{ to } j .&\end{array}\right. 
\end{equation}

With this in place, the second main result of the paper can be formalized. 

\noindent\colorbox{mygray}{\begin{minipage}{0.48\textwidth}\begin{theorem}\label{policy}
Let each computing node $i$ use the policy $\pi_{ij}$ (for each $j$, such that $(i,j)\in E$) defined in~\eqref{pi} for specified values of $\epsilon, \eta$, and $\nu$. Under Assumptions~\ref{as.function} till Assumption~\ref{as.persistent} and with the same definitions of Theorem~\ref{convergence}, the quantity $\tilde{\g}_{i,k}$ defined in~\eqref{gtilde} satisfies $\|\tilde{\g}_{i,k} - \g_{i,k}\| \leq \epsilon/2 |X_{k,k-1}|$ and it is a stochastic $\epsilon$-(sub)gradient of $f_{i,k}(\x,\myom)$ w.r.t. $\v_{i,k+1}$. 
\end{theorem}\end{minipage}}

\vskip0.1cm Theorem~\ref{policy} guarantees that our proposed utility-based policy will deliver $\epsilon$-gradient vectors of a prescribed accuracy, thereby ensuring convergence in the sense of Theorem~\ref{convergence}. The policy $\pi_{ij}$ is completely disconnected with the running of Algorithm~\ref{alg}, and it can be triggered at any time step $k$. Its utility-based nature allows each node to send its neighbors the up-to-date information only when they really need it. We remark that the information that each node has to send is the gradient of $f_{i,k,j}$ with respect to $\x$ at $\v_{i,k+1}$, and the \textsc{pdf} of $\myom_{i,k}$. Both are functions of $\myom_{i,k}$, which need to be encoded in some way. We will show with a simple numerical evaluation how this can be done in practice and that the added communication cost in sending $\nabla_{\x} f_{i,k,j}(\x,\myom_j)|_{\v_{i,k+1}}$ (which is not present if we send the \textsc{pdf} of $\myom_{i,k}$ at each time $k$) is completely justifiable.  

\begin{remark}
The need of sending $\nabla_{\x} f_{i,k,j}(\x,\myom_j)|_{\v_{i,k+1}}$ can be properly tuned by selecting small or high $\nu$. In the limit $\nu \to 0$, there is no need of sending it. 
\end{remark}

\begin{remark}
In the policy $\pi_{ij}$, each node needs to know $\Delta_i$ and $\Delta_j$ for all its neighboring nodes. This can be avoided by considering a more conservative policy which substitutes $\Delta_i$ and $\Delta_j$ with $n$. 
\end{remark}

\begin{remark}\label{linear}
If each of the functions $f_{i,k}(\x,\myom)$ is a sum of two functions, one dependent on $\x$ and one \emph{linearly} dependent on $\myom$, then it is easy to see that the policy $\pi_{ij}$ becomes checking the utility metric $\mathcal{U}_{\mathrm{S},ij}^{(1)}(\myom_k|\myom')$ only, as we will see in one of the numerical examples.     
\end{remark}

\section{Numerical Evaluations}\label{sec:num}

We now look at two different numerical examples to assess the proposed algorithm in terms of communication overhead and general performance.  %The aim is to highlight the main properties and gain insights on the actual savings in terms of number of messages. 

\subsection{Least-squares estimation in sensor networks}

As a first example, we consider an estimation problem, where each of the computing nodes needs to estimate the state $\x\in[-\nicefrac{1}{2},\nicefrac{1}{2}]^d$ of a slowly time-varying process. 
%a linear discrete time-invariant process,
%\begin{equation*}
%\x_{k+1} = {\bf \Phi}\x_k + {\bf n},
%\end{equation*}
%where ${\bf \Phi}$ is the state transition matrix, while ${\bf n}$ is the process noise, assumed zero-mean gaussian with covariance ${\bf Q}$. 
Each node relies on its own measurements ${\bf z}_{i,k} \in\mathbb{R}^d$ that are linked to the state via the measurement equation, 
\begin{equation}\label{meas}
{\bf z}_{i,k} = {\bf H}_i(\myom)\x_k + {\bf w}_i,
\end{equation}
where ${\bf H}_i$ is a squared stochastic matrix assumed to be full rank for the sake of simplicity. The noise term ${\bf w}_i$ is also zero-mean Gaussian and has covariance ${\bf R}_i$. We assume that each of the ${\bf H}_i$'s has the form,
\begin{equation}\label{coupled}
{\bf H}_i(\myom) = \bar{\bf H}_i + {\bf{I}}_d \sum_{j\in\myset{N}_i\cup i} c_j \omega_j,
\end{equation}
while we assume that each scalar $\omega_i\in\mathbb{R}_{+}$ is drawn from a truncated Rayleigh distribution with scale parameter $\sigma_{i,k} > 0$ over the compact set $[0, 3]$, slowly varying with 
time according to a truncated first-order model, 
$$
\sigma_{i,k+1} = \proj_{[0,3]}\left[\sigma_{i,k} + \rho_i \sin(a_i k + b_i) + r_i\right], \quad r_i \sim \mathcal{N}(0, P_i).  
$$

We consider the following time-varying least-squares estimator for the process $\x$ at each time $k$
\begin{equation}\label{est}
\minimize_{\x_k \in [-\nicefrac{1}{2},\nicefrac{1}{2}]^d} \myexp_{\myom,k} \Big[\sum_{i=1}^n\|{\bf z}_{i,k} - {\bf H}_i(\myom)\x_{k}\|^2_{{\bf R}_i^{-1}}\Big]
\end{equation}
and we approximate the expectation operator as a finite sum over $N$ Monte Carlo samples of the stochastic variables $\myom$.

Problem~\eqref{est} is the optimization program that the computing nodes have to solve in a distributed way using the proposed algorithm, and it is indeed a particular instance of~\eqref{optimalproblem}. Distributed estimation tasks of the type of~\eqref{est} may arise in synthetic aperture radar (SAR) array scenarios, for example~\cite{Pardini2008, Bioucas-Dias2010}. In SAR arrays, noise can enter in the measurements in a multiplicative fashion (speckle noise) and it can be modeled with a Rayleigh \textsc{pdf}. A coupled measurement equation with of the form~\eqref{meas}-\eqref{coupled} can be used to model correlated noise within the neighboring sensors.  Other similar problems arise when sensors have access to pair-wise relative measurements.

The simulation that we are going to present has the following parameters\footnote{The code of the simulation will be made available on-line, so that the interested reader can reproduce the results and change the parameters at will.}: the dimension of the state vector is $d=2$, the number of sensors is $n = 15$, the true state evolves via a linear discrete time-invariant process (unknown to the sensor nodes),
\begin{equation*}
\x^{\textrm{true}}_{k+1} = {\bf \Phi}\x^{\textrm{true}}_{k} + {\bf n},
\end{equation*}
where ${\bf \Phi}$ is the state transition matrix, while ${\bf n}$ is the process noise, assumed zero-mean gaussian with covariance ${\bf Q}$ and in particular, ${\bf \Phi} = [.99, 0.01; 0, 1]$, ${\bf Q} = \xi_Q{\bf I}_2$,  $\bar{\bf H}_i = {\bf I}_2$, ${\bf R}_i = \xi_R{\bf I}_2$, $\xi_Q = \xi_R = 1$e-$6$, $\sigma_{i,0} = \max\{i/n + 0.3(\mathcal{U}_{[0,1]}-0.3),0.001\}$, $P_i = 1$e-$2$, $a_i = i \pi /200$, $b_i= 200\pi\rho_i$, $\rho_i \sim \mathcal{N}(0,1)$, and the number of Monte-Carlo samples for the stochastic variable $\myom$ is $N=5000$. We also consider $s_{ij}$ to be the same for all the links. Also, we scale~\eqref{est} by multiplying by $\xi_R$ for better numerical stability. For the problem at hand, we can easily compute the bound on $\beta$ (we set $\beta = 1/n-1$e-$4$), while a suitable value for $\alpha$ has to be derived by trial-and-error as usually happens in distributed optimization (since the constants $m_f$ and $L$ are difficult to obtain in practice). Guided by the stochastic-free case, we set $\alpha= 1/400$ which works well in the simulation scenarios. %\footnote{Given the nature of the cost function and the fact that the stochastic variables lie in a compact set, a conservative bound can be also derived analytically. However, in this particular example, this bound turns out to be extremely conservative and it has not been used to determine $\alpha$.}. 
%and
%$$
%\beta \leq \nicefrac{s_{ij}}{(s_{ij} \max_i\{\Delta_i\} + 1)}
%$$
%
Furthermore, we can explicitly write the gradient $\nabla_{\x} f_{i,k}(\x,\omega_j)|_{\v_{i,k+1}}$ of the scaled problem as, 
\begin{equation*}
\nabla_{\x} f_{i,k}(\x,\omega_j)|_{\v_{i,k+1}} = 2\left( {\bf{H}}_i(\myom)^\transp {\bf{H}}_i(\myom)\v_{i,k+1} - {\bf z}_{i,k}^\transp {\bf{H}}_i(\myom)\right)
\end{equation*}
so that with the policy $\pi_{ij}$ we need to send either the scalar parameter $\sigma_{i,k}$ or the couple $(\v_{i,k+1},{\bf z}_{i,k})$, thus $4$ scalar values.  Finally, the optimal solution of~\eqref{est} for comparison is also computable in closed form. 

%\textbf{Example 1.} In the first example, 
In this simulation example, we look at loosely connected, highly asynchronous sensor nodes with $s_{ij} = 0.3$, $\gamma = 0.98$, and a high noise dependence with $c_i = 1$ in ${\bf H}_i$ for all sensors. 

\begin{figure}
\centering
\psfrag{a}{\footnotesize \hspace{-1cm} discrete time $k$}
\psfrag{b}{\footnotesize \hspace{-1cm} $\myexp[\|\y_k - \y^*_k\|^2]$}
\psfrag{c}{\footnotesize $\epsilon = 0.001, \nu = .25\epsilon$}
\psfrag{d}{\footnotesize $\epsilon = 5, \nu = .25\epsilon$}
\psfrag{e}{\footnotesize No \textsc{pdf} communication}
\psfrag{x}{\footnotesize \hspace{-2.5cm} Est. bound, $\epsilon=0.001$: Theorem~\ref{convergence}}
%\psfrag{d}{\footnotesize $\epsilon = \nu = 0.001$}
%\psfrag{e}{\footnotesize $\epsilon = \nu = 1$}
%\psfrag{f}{\footnotesize No communication}
%\psfrag{g}{\footnotesize Error bound for $\epsilon = \nu = 0.001$ }
\includegraphics[width=0.50\textwidth, trim=.75cm 0.5cm .75cm .5cm, clip=on]{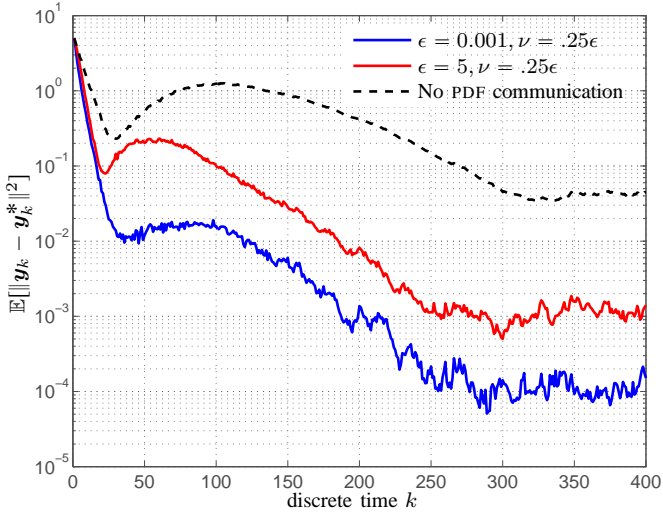} 
\caption{Distributed least-squares: performance in terms of error w.r.t. the optimizer. The results have been averaged over $25$ realizations.}\label{fig1a}
\end{figure}

\begin{figure}
\centering
\psfrag{a}{\footnotesize \hspace{-1cm} discrete time $k$}
%\psfrag{b}{\footnotesize \hspace{-1cm} $\myexp[\|\y_k - \y^*_k\|^2]$}
\psfrag{b}{\footnotesize \hspace{-2.3cm} Utility based/ every-time communication}
\psfrag{c}{\footnotesize $\epsilon = 0.001, \nu = .25\epsilon$}
\psfrag{d}{\footnotesize $\epsilon = 5, \nu = .25\epsilon$}
\psfrag{e}{\footnotesize No \textsc{pdf} communication}
%\psfrag{g}{\footnotesize Error bound for $\epsilon = \nu = 0.001$ }
\includegraphics[width=0.50\textwidth, trim=.75cm 0.5cm .75cm .5cm, clip=on]{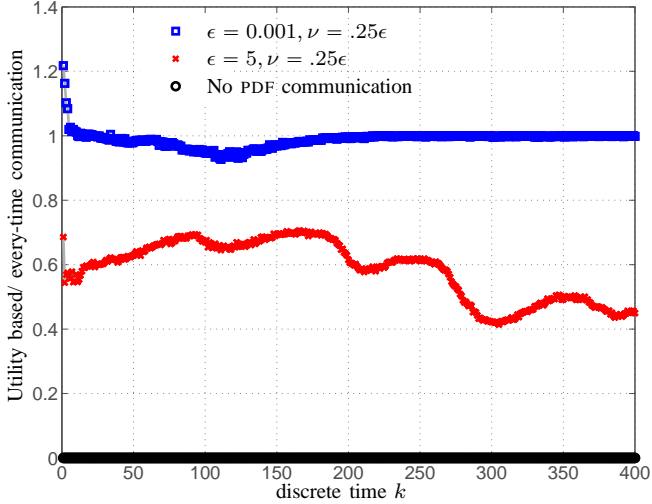}
\caption{Distributed least-squares: performance in terms of \textsc{pdf} communication. The results have been averaged over $25$ realizations.}\label{fig1b}
\end{figure}

The results show convergence of the proposed scheme to an error floor dependent on $\epsilon$ (Figure~\ref{fig1a}). %We have also displayed the theoretical bound of Theorem~\ref{convergence} (for $\epsilon=0.001$), for which $L$ and $m_f$ are estimated a posteriori on the data. This shows that the bound is reasonably accurate for small $\epsilon$ and $\delta_\x$. 
Furthermore, it is possible to appreciate the decrease in communication overhead when selecting different error levels w.r.t. the every-time communication strategy, where each node sends to all its neighbors its updated \textsc{pdf} at every time $k$  (Figure~\ref{fig1b}). In particular, even selecting $\epsilon=5, \nu = .25\epsilon$ yields acceptable accuracy with a significant reduction of communication.  (We notice that the utility-based communication can have a communication count greater than every-time, since in utility-based the nodes may need to send also the gradient of the local cost
 functions).

%\textbf{Example 2.} The second example deals with more connected nodes, more synchronous communication, $s_{ij} = 0.7$, $\gamma = 0.32$, and low noise dependence with $c_i = 0.5/n$. In this case, we see that the dependence on the stochastic variable is less important than in the first example, and as a consequence, we can obtain good accuracies even with higher error levels (Figure~\ref{fig2}).

\begin{remark}\label{remdin1}
More advanced examples can be derived from this least-squares estimation problem. For instance, moving horizon strategies can be considered in the context of maximum a posteriori estimation, as in~\cite{Jakubiec2013}. One thing to pay attention to is the possible dependence of the cost function on past (estimated) values of the state, i.e., $\y_{i,k-1}$; in this case, although the proposed algorithm still converges (under the given assumptions), the physical meaning of the limit point may be less clear. 
\end{remark}

%\begin{figure}[h]
%\centering
%\psfrag{a}{\footnotesize \hspace{-1cm} discrete time $k$}
%\psfrag{b}{\footnotesize \hspace{-1cm} $\myexp[\|\y_k - \y^*_k\|^2]$}
%\psfrag{c}{\footnotesize \hspace{-2.3cm} Utility based/ every-time communication}
%\psfrag{d}{\footnotesize $\epsilon = \nu = 0.001$}
%\psfrag{e}{\footnotesize $\epsilon = \nu = 0.1$}
%\psfrag{f}{\footnotesize No communication}
%\psfrag{g}{\footnotesize Error bound for $\epsilon = \nu = 0.001$ }
%\includegraphics[width=0.50\textwidth]{figure/g2error} \\ \includegraphics[width=0.50\textwidth]{figure/g2comm}
%\caption{\textbf{Example 2.} Highly connected scenarios for different values of $\epsilon$, $s_{ij} = 0.7$, $\gamma = 0.32$, $c_i = 0.5/n$. The results have been averaged over $50$ realizations.}\label{fig2}
%\end{figure}

\subsection{Waypoint generation in robotic networks}

The second example is a networked formation control problem, where a number of mobile nodes needs to track a defined point in space and maintain a certain formation. The example is inspired by~\cite{Borrelli2008} and has the added aim to show (\emph{i}) that the proposed algorithm can work with partially overlapping decision variables $\x$ (i.e., there is no need for each of the computing nodes to agree on the total decision variable $\x$ but only on subsets of it), and (\emph{ii}) how the policy $\pi_{ij}$ gets simplified in the case of noise entering linearly in the cost function (see Remark
~\ref{linear}). 

\begin{figure}
\centering
\psfrag{a}{\footnotesize \hspace{-0.5cm} $k = 2$}
\psfrag{b}{\footnotesize \hspace{-0.5cm} $k = 100$}
\psfrag{c}{\footnotesize \hspace{-0.5cm} $k = 250$}
\psfrag{d}{\footnotesize \hspace{-0.6cm} $k = 500$}
\psfrag{e}{\footnotesize \hspace{-0.6cm} $k = 1000$}
\psfrag{f}{\footnotesize \hspace{-0.6cm} $k = 5000$}
%\psfrag{g}{\footnotesize Error bound for $\epsilon = \nu = 0.001$ }
\noindent\includegraphics[width=0.50\textwidth, trim=2cm .5cm 1.25cm .5cm, clip=on]{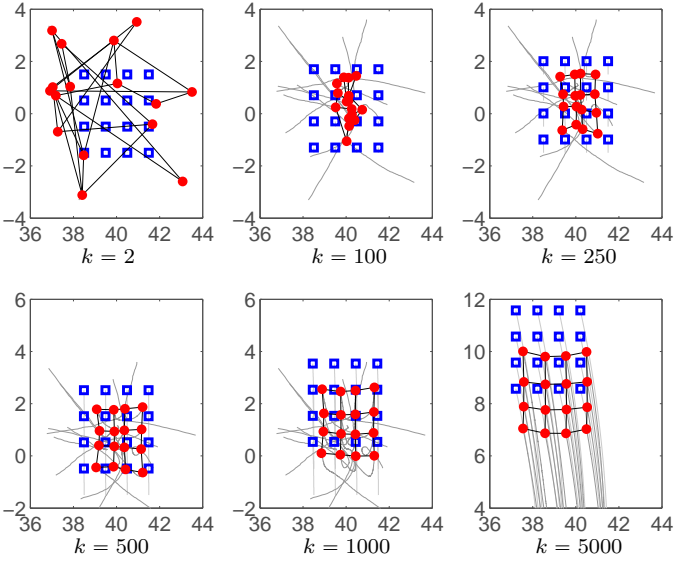}
\caption{Snapshots of the algorithm's waypoint generation (red points) and reference ones (blue squares) in the case of no noise. }\label{fig3}
\end{figure}

We consider $n = 16$ mobile nodes that have a fixed connection structure and need to track a squared pattern figure in two dimensions (Figure~\ref{fig3}). At a given discrete time $k$, each mobile node $i$ needs to compute a waypoint $\x_{i,k}$ where to head to, this waypoint depends on the current value of the reference point ${\bf x}^\mathrm{ref}_{i,k}$ and on the neighboring waypoint/reference values. In addition, we consider each of the reference points to be known to the nodes with some degree of uncertainty, as 

$$
{\x}^\mathrm{ref}_{i,k} = \bar{{\x}}^\mathrm{ref}_{i,k} + \myom_{i,k},
$$ 
where $\bar{{\bf x}}^\mathrm{ref}_{i,k}$ is the actual reference and $\myom_{i,k}$ is drawn from a given (estimated) \textsc{pdf}. Putting this together, the computing mobile nodes have to solve the optimization problem
\begin{multline}\label{control}
\minimize_{{\x}_{k} \in X_k}\mathbb{E}_{\myom,k} \Big[\sum_{i\in V} \Big(\theta \|{\x}_{i,k}-{\x}^\mathrm{ref}_{i,k}\|^2 + \\
\sum_{j\in N_i}\|{\x}_{i,k} - {\x}_{j,k}-({\x}^\mathrm{ref}_{i,k}-{\x}^\mathrm{ref}_{j,k})\|^2\Big)\Big],
\end{multline}
where $\x_{k}$ is the stacked version of the all $\x_{i,k}$'s and $\theta>0$ a chosen scaling factor.  
%
%
%
%The position of each node is denoted as ${\bf p}_{i,k}\in\mathbb{R}^2$, and evolves via the first-order discrete-time linear system
%$$
%{\bf p}_{i,k+1} = {\bf p}_{i,k} + {\bf u}_{i,k} + \myom_{i,k}, 
%$$ 
%where ${\bf u}_{i,k}\in U_{i,k}\subseteq \mathbb{R}^2$ is the control action and $\myom_{i,k}$ is a disturbance. We assume that the disturbance has a fixed direction in space but can vary its intensity, that is $\myom_{i,k} = (\omega_{i,k}; \omega_{i,k})$. The aim of the network of nodes is to solve the following optimization program at each discrete time step
%\begin{multline}\label{control}
%\minimize_{{\bf u}_{k} \in U_{k}}\mathbb{E}_{\myom,k} \Big[\sum_{i\in V} \|{\bf p}_{i,k+1}-{\bf p}^\mathrm{ref}_{i,k+1}\|^2_{{\bf Q}^{-1}} + \\
%\sum_{(i,j)\in E}\|{\bf p}_{i,k+1} - {\bf p}_{j,k+1}-({\bf p}^\mathrm{ref}_{i,k+1}-{\bf p}^\mathrm{ref}_{j,k+1})\|^2_{{\bf R}^{-1}}\Big],
%\end{multline}
%where.. 
%

Problem~\eqref{control} is an instance of~\eqref{optimalproblem}.  Problems like~\eqref{control} arise in multi-agent scenarios where each mobile node reference point knowledge is affected by stochastic disturbances. This added term can be a pre-compensation (before running the trajectory generation algorithm) of known disturbances acting during the time frame in which the mobile node aims at reaching the waypoint. These disturbances can be internal (e.g., motors), or external (e.g., environment). Environment noise can be for example wind, if the nodes are sufficiently far apart so that the wind variations at the nodes w.r.t. the average value are uncorrelated with each other. Wind is often modeled as a Weibull distribution, and therefore with this in mind, we model each of the stochastic variables $\myom_{i,k}$'s as drawn from a  Weibull \textsc{pdf} varying in time (there is no need for the compactness assumption for $\Omega$ in this case, as we will see). In particular, $\myom_{i,k} = {\bf 1}_2 \omega_{i,k}$, $\omega_{i,k} \sim \mathcal{W}(\lambda_{i,k}, \sigma_{i,k})$,   
$$
\lambda_{i,k} = \lambda_{i,k-1}(1 + \cos(r + \phi_{\lambda,i})/1.3) + \omega_{\max,i}, 
$$
$$
\sigma_{i,k} = \sigma_{i,k-1}(1 + \cos(r + \phi_{\sigma,i})/1.3) + \omega_{\max,i}, 
$$
where, $\mathcal{W}(\lambda, \sigma)$ is the Weibull distribution with scalar parameters $\lambda$ and $\sigma$, $\phi_{\lambda,i}$ and $\phi_{\sigma,i}$ are drawn from $\mathcal{U}_{[0, 2\pi]}$, and $r$ and $\omega_{\max,i}$ are given scalars.  

The reference states ($\bar{\x}_{i,k}^{\mathrm{ref}}$) evolve along circular trajectories with constant angular velocity $r$. The other parameters of the simulation example are: $s_{ij} = 0.7$, $\theta = 0.5$, $X_k = [-50,50]^{32}$ $\sigma_{i,0} = 4 + \mathcal{U}_{[0,4]}$, $\lambda_{i,0} = 4/n \omega_{i,\max}(2 + \mathcal{U}_{[0,1]})$. The step sizes $\alpha$ and $\beta$ are determined according to Theorem~\ref{convergence}, in fact, in this example, the bounds are analytically computable. For the simulations we set $\alpha = 1.7 \theta/(2(\theta + \lambda_n(\bar{\W})/s_{ij}))^2$ and $\beta = 0.13/n$. We can write the gradient $\nabla f_{i,k}$ in closed form as 
$$
\nabla_{\x} f_{i,k} = 2 \theta\big({\bf x}_{i,k}-{\bf x}^\mathrm{ref}_{i,k}\big) +
2 \sum_{j\in N_i}{\bf x}_{i,k} - {\bf x}_{j,k}-({\bf x}^\mathrm{ref}_{i,k}-{\bf x}^\mathrm{ref}_{j,k}),
$$ 
and since the gradient depends linearly on $\myom_{k}$ (through ${\x}_{i,k}^{\mathrm{ref}}$), we can simplify the policy $\pi_{ij}$ by noticing that
\begin{equation}\label{simplcond}
\|\tilde{\g}_{i,k} -\g_{i,k}\| \leq \sum_{j\in N_i} 2\|\widehat{\myom_{j}'}-\widehat{\myom_{j,k}}\| \leq \epsilon/2|X_{k,k-1}|,
\end{equation}
where $\widehat{\cdot}$ represents the mean operator; therefore the policy $\pi_{ij}$ boils down to sending the updated \textsc{pdf} only when condition~\eqref{simplcond} does not hold.  

\begin{figure}
\centering
\psfrag{a}{\footnotesize Bound: Theorem~\ref{convergence}}
%\psfrag{b}{\footnotesize $\alpha'' = \alpha, \beta'' = 10\beta$}
%\psfrag{c}{\footnotesize  $\alpha''' = \alpha/10, \beta''' = \beta$}
\psfrag{x}{\footnotesize \hspace{-1cm} discrete time $k$}
\psfrag{y}{\footnotesize \hspace{-1cm} $\myexp[\|\y_k - \y^*_k\|^2]$}
%\psfrag{g}{\footnotesize Error bound for $\epsilon = \nu = 0.001$ }
\noindent\includegraphics[width=0.50\textwidth, trim=.75cm 0.5cm .75cm .0cm, clip=on]{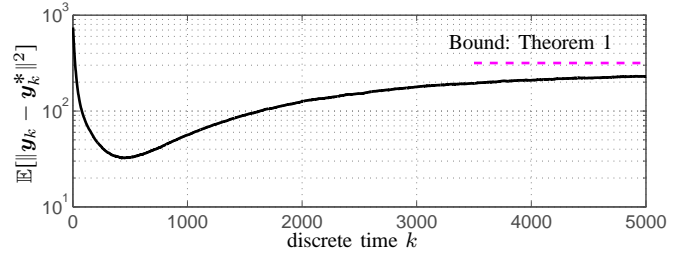}
\caption{Convergence performance of the algorithm without noise.}\label{fig4}
\end{figure}

\textbf{Example 1.} In the first example we consider a simulation scenario without noise. The purpose is to display the performance of the algorithm in its ideal case, when both the gradient is computed exactly and the stochastic program reduces to a deterministic optimization problem (and no compactness assumption on $X_k$ is required). We select the angular velocity as $r = 0.5/40 \alpha$, and run the distributed asynchronous time-varying optimization problem up to $k = 5000$. By using snapshots of the agents trajectories we show the algorithm's behavior (Figure~\ref{fig3}). The blue squares are the reference waypoints, while the red points are the agent-computed waypoints at discrete time $k$. The lines represent the possible connection among agents (which are not time-varying~\cite{Borrelli2008}). As we further see, the convergence performance is in line with the asymptotical bound of Theorem~\ref{convergence}, which is rather tight in this particular case (Figure~\ref{fig4}).

\begin{figure}
\centering
\psfrag{a}{\footnotesize \hspace{-1cm} discrete time $k$}
\psfrag{b}{\footnotesize \hspace{-1cm} $\myexp[\|\y_k - \y^*_k\|^2]$}
\psfrag{c}{\footnotesize $\epsilon = 0.02$}
\psfrag{d}{\footnotesize $\epsilon = 0.2$}
\psfrag{e}{\footnotesize No \textsc{pdf} communication}
%\psfrag{}{\footnotesize \hspace{-1cm} Approx. Bound $\epsilon = 1$: Lemma~\ref{lemma.3}}
\psfrag{x}{\footnotesize \hspace{-.5cm} Bound $\epsilon = .02$ $\uparrow$}
\psfrag{y}{\footnotesize \hspace{-.5cm} Bound $\epsilon = .2$ $\downarrow$}
%\psfrag{g}{\footnotesize Error bound for $\epsilon = \nu = 0.001$ }
\noindent\includegraphics[width=0.50\textwidth, trim=.75cm 0.5cm .75cm .5cm, clip=on]{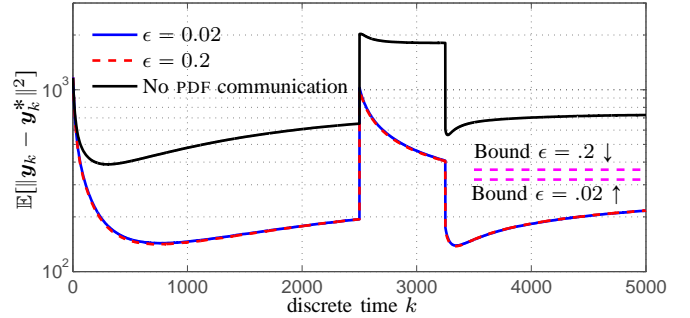}
\caption{Distributed waypoint generation: performance in terms of error w.r.t. the optimizer. The results have been averaged over $25$ realizations.}\label{fig5}
\end{figure}

\begin{figure}
\centering
\psfrag{x}{\footnotesize \hspace{-1cm} discrete time $k$}
\psfrag{y}{\footnotesize \hspace{-2.3cm} Utility based/ every-time communication}
\psfrag{c}{\footnotesize $\epsilon = 0.02$}
\psfrag{d}{\footnotesize $\epsilon = 0.2$}
\psfrag{e}{\footnotesize No \textsc{pdf} communication}
%\psfrag{g}{\footnotesize Error bound for $\epsilon = \nu = 0.001$ }
\noindent\includegraphics[width=0.50\textwidth, trim=.75cm 0.5cm .75cm .5cm, clip=on]{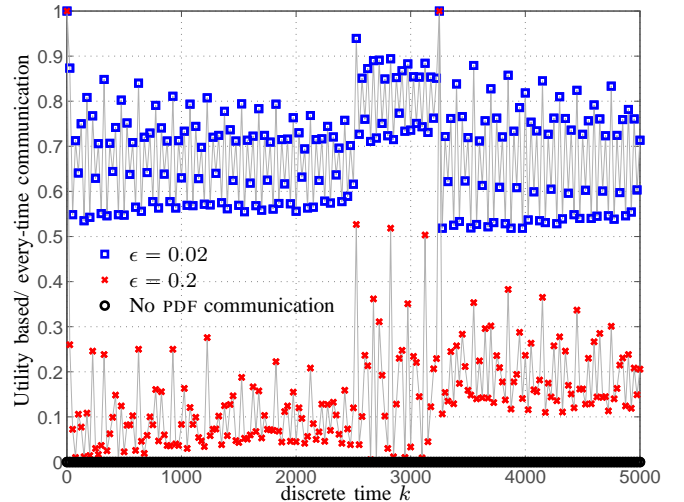}
\caption{Distributed waypoint generation: performance in terms of \textsc{pdf} communication. The results have been averaged over $25$ realizations. (Only one out of every $25$ time instances is shown for clarity).}  \label{fig6}
\end{figure}

\textbf{Example 2.} Stochastic noise is introduced in the second simulation example. We consider the same angular velocity and same $\alpha$ and $\beta$. We select $\omega_{i,\max}$ as $2 i /n $ and we run Algorithm~\ref{alg} for different choices of $\epsilon$ for the utility-based policy $\pi_{ij}$. To simulate a sudden change in the noise (e.g., a wind gust), we double $\lambda_{i,k}$ in the range $k\in[2500, 3250]$. As we can see in Figures~\ref{fig5} and \ref{fig6}, even considering limited
 communication exchange,  Algorithm~1 with policy $\pi_{ij}$ performs almost as if all the nodes were communicating their \textsc{pdf} at each time step. We see that the bound of Theorem~\ref{convergence} for $\epsilon=.02$ is tighter than for $\epsilon=.2$, as expected (see Remark~\ref{tight}); nonetheless, the utility-based policy works as designed.

%\clearpage
%\newpage

\section{Proofs} \label{proofs}

\subsection{Setting Up the Analysis}

We define $\g_{i,k}$ as the stochastic gradient of $f_{i,k}$ at $\v_{i,k+1}$,
\begin{equation}\label{gom}
\g_{i,k} = \myexp_{\myom,k}[\nabla_\x f_{i,k}(\x,\myom)|_{\v_{i,k+1}}],
\end{equation}
a quantity that will be used to characterize how large ``$\epsilon$'' is for the stochastic $\epsilon$-(sub)gradient $\tilde{\g}_{i,k}$. We also define the stochastic gradient of $f_{i,k}$ at the optimizer of~\eqref{optimalproblem}, i.e., $\x^{*}_{k} = \y_{i,k}^* = \v_{i,k+1}^*$, as,
\begin{equation}\label{gopt}
\g_{i,k}^* = \myexp_{\myom,k}[\nabla_\x f_{i,k}(\x,\myom)|_{\y^*_{i,k}}].
\end{equation}

\noindent\colorbox{mygray}{\begin{minipage}{0.48\textwidth}\begin{lemma}\label{lemma1}
Let $X_{k,k-1}$ be defined as in Theorem~\ref{convergence}; let $\v_{i,k+1}$ be defined as \eqref{v}, let $\tilde{\g}_{i,k}$, $\g_{i,k}$ and ${\g}_{i,k}^*$ be defined as~\eqref{gtilde}, \eqref{gom}, and \eqref{gopt}, respectively. Under Assumption~\ref{as.function} and for $\beta < 1/n$, then 
\begin{enumerate}
\item[\emph{(a)}] the supporting variable $\v_{i,k+1} \in \myset{X}_{k-1}$;
\item[\emph{(b)}] if $\|\tilde{\g}_{i,k}-\g_{i,k}\| \leq \epsilon/2 |X_{k,k-1}|$ holds, then $\tilde{\g}_{i,k}$ is a stochastic $\epsilon$-(sub)gradient of $f_{i,k}(\x,\myom)$ at $\v_{i,k+1}$; 
\item[\emph{(c)}] at each discrete time step $k$, the inequality $\|{\g}_{i,k}-\g_{i,k}^*\| \leq L\|\v_{i,k+1} - \y^{*}_{i,k}\|$ holds.
\end{enumerate}
\end{lemma}\end{minipage}}
\vskip0.1cm \noindent \begin{myproof} 
\emph{(a)} The variable $\v_{i,k+1}$ is generated in~\eqref{v} as
\begin{equation*}
\v_{i,k+1} \hskip-0.1cm= \hskip-0.1cm%\y_{i,k} - \beta \sum_{j = 1}^n [\W_k]_{i,j} \y_{j,k} =\\= 
\Big(1 - \beta\hskip-0.3cm \sum_{j=1, j\neq i}^n\hskip-0.2cm [{\bf A}_k]_{i,j} \Big)\y_{i,k} + \beta\hskip-0.3cm\sum_{j=1, j\neq i}^n \hskip-0.2cm[{\bf A}_k]_{i,j} \y_{j,k} \hskip-0.1cm = \hskip-0.1cm\sum_{j=1}^n \theta_j \y_{j,k},
\end{equation*}
with $\sum_{j=1}^n \theta_j = 1$, $\theta_{j,j\neq i}\geq 0$. Furthermore, %for the weighted largest eigenvalue of Laplacian graph holds the lower bound (coming directly from its definition)~\cite{Berman2011}, 
%$$
%\bar{s} := \max_{i\in\myset{V}}\Big\{\sum_{j=1, j\neq i}^n s_{ij} \Big\} \leq \lambda_n(\bar{\W}).
%$$
%This implies for each $i$
$$
\sum_{j=1, j\neq i}^n [{\bf A}_k]_{i,j} < n < 1/\beta,
$$
and thus it is also true that $\theta_i  = 1 - \beta \sum_{j=1, j\neq i}^n [{\bf A}_k]_{i,j} > 0$. This means that $\v_{i,k+1}$ is a convex combination of $\y_{i,k}$ and since $\y_{i,k}\in\myset{X}_{k-1}$, then $\v_{i,k+1}\in\myset{X}_{k-1}$.

\emph{(b)} By the definition of the $\epsilon$-(sub)gradient given in~\eqref{egrad} and by the definition of the gradient, by the fact that the domain of the $f_{i,k}$'s needs to be at least $X_{k,k-1
}\times\Omega$ by assumption, it has to be
\begin{multline}\label{dum}
(\tilde{\g}_{i,k} - \g_{i,k})^\transp (\v-\v_{i,k+1}) + \g_{i,k}^\transp (\v-\v_{i,k+1})\leq \\ \myexp_{\myom,k}[f_{i,k}(\v,\myom) - f_{i,k}(\v_{i,k+1},\myom)]+ \epsilon,
\end{multline}
for all ${\v} \in X_{k,k-1}$, $\epsilon \geq 0$. Since by definition, the gradient verifies the same condition of~\eqref{egrad} with $\epsilon= 0$, a sufficient condition for~\eqref{dum} is
$$
\left\|(\tilde{\g}_{i,k} - \g_{i,k})^\transp (\v-\v_{i,k+1})\right\| \leq \epsilon,
$$
and yet another sufficient condition for this latter inequality to hold is
$$
\|\tilde{\g}_{i,k} - \g_{i,k}\| \leq \frac{\epsilon}{2 |X_{k,k-1}|}.
$$

\emph{(c)} This claimed inequality comes from the Lipschitz assumption applied to the gradients $\g_{i,k}$ and $\g_{i,k}^*$ of the points $\v_{i,k+1}\in\myset{X}_{k-1}$ and $\y_{i,k}^* \in\myset{X}_{k-1}$. 
\end{myproof}

\noindent\colorbox{mygray}{\begin{minipage}{0.48\textwidth}\begin{lemma}\label{lemma.2}
With the same definitions of Theorem~\ref{convergence} and under Assumption~\ref{as.function}, for any $\mu >0$, if $\|\tilde{\g}_{i,k}-\g_{i,k}\| \leq \epsilon/2 |X_{k,k-1}|$ for all $i,k\geq 1$, then the sequence $\{\y_{i,k}\}$ generated by the proposed algorithm is bounded as 
\begin{multline*}
\hskip-0.35cm\|\y_{i,k+1} - \y_{i,k}^{*}\|^2 \hskip-0.1cm\leq \hskip-0.1cm(1+\mu)(1+\alpha^2L^2- \alpha m_f)\left\|\v_{i,k+1}-  \y_{i,k}^{*}\right\|^2 + \\
 (1+\mu) \alpha \epsilon \left(\alpha\frac{\epsilon}{4 |X_{k,k-1}|^2} + 2 \alpha L  + 2 \right) + (1+1/\mu)\alpha^2 G^2.
\end{multline*}
\end{lemma}\end{minipage}}\vskip0.1cm
%------------------------------------------------------------------------------------
%------------------------------------------------------------------------------------
%------------------------------------------------------------------------------------
\noindent\begin{myproof} 
By the definition of the update rule~\eqref{update} and the optimal gradient $\g_{i,k}^*$ given in~\eqref{gopt}, we can write 
\begin{multline}\label{lemma.eq1}
\hskip-0.3cm\|\y_{i,k+1} - \y_{i, k}^{*}\|^2 \hskip-0.1cm= \hskip-0.1cm \big\|\proj_{\myset{X}_{k}}[\v_{i,k+1} - \alpha \tilde{\g}_{i,k}] -\proj_{\myset{X}_{k}} [\y_{i,k}^{*}-\alpha \sum_{j\in\myset{V}}{\g}_{j,k}^*]\big\|^2 \\
\leq \Big\|\v_{i,k+1} - \alpha (\tilde{\g}_{i,k}-\g_{i,k}^*)-  \y_{i,k}^{*} + \alpha\sum_{j\in\myset{V}, j\neq i}{\g}_{j,k}^*\Big\|^2,
\end{multline}
%--
where we have used $\y_{i,k}^{*} = \proj_{\myset{X}_{k}} [\y_{i,k}^{*}-\alpha \sum_{j\in\myset{V}}{\g}_{j,k}^*]$, and the non-expansivity property of the projection. If we use the Peter-Paul inequality\footnote{I.e., $(a+b)^2\leq (1+\mu)a^2 + (1+1/\mu)b^2$, for all $a,b \in\mathbb{R}$ and $\mu>0$.} and the boundedness of the optimal gradient~(Eq.~\eqref{b.G}), we obtain
\begin{multline}\label{lemma.eq1bis}
\hskip-0.3cm\|\y_{i,k+1} - \y_{i, k}^{*}\|^2 \leq (1+\mu)\|\v_{i,k+1} - \alpha (\tilde{\g}_{i,k}-\g_{i,k}^*)-  \y_{i,k}^{*}\|^2 +\\ (1+1/\mu)\alpha^2 G^2,
\end{multline}
If we expand the first part of the bound in~\eqref{lemma.eq1bis}, we obtain
\begin{multline}\label{lemma.eq2}
%\|\y_{i,k+1} - \y_{i,k}^{*}\|^2 \leq 
\left\|\v_{i,k+1} - \alpha (\tilde{\g}_{i,k}-\g_{i,k}^*)  -  \y_{i,k}^{*}\right\|^2 = 
\left\|\v_{i,k+1}-  \y_{i,k}^{*}\right\|^2 + \\ \alpha^2\left\|\tilde{\g}_{i,k}-\g_{i,k}^*\right\|^2 - 2 \alpha(\tilde{\g}_{i,k}-\g_{i,k}^*)^\transp (\v_{i,k+1}-  \y_{i,k}^{*}). 
\end{multline}
By using Lemma~\ref{lemma1}, Assumption~\ref{as.function}, and the definition of $\epsilon$-(sub)gradient, we can now bound the following terms, 
\begin{multline*}
\left\|\tilde{\g}_{i,k}-\g_{i,k}^*\right\|^2 \leq (\|\tilde{\g}_{i,k}-\g_{i,k}\| + \|{\g}_{i,k}-\g_{i,k}^*\|)^2\leq \\ \Big(\frac{\epsilon}{2|X_{k,k-1}|} + L \|\v_{i,k+1}-\y_{i,k}^*\|\Big)^2 \leq \\ \frac{\epsilon^2}{4 |X_{k,k-1}|^2} + 2 L \epsilon + L^2 \|\v_{i,k+1}-\y_{i,k}^*\|^2
\end{multline*}
where we bounded $\|\v_{i,k+1}-\y_{i,k}^*\|\leq 2 |X_{k,k-1}|$, and by strong convexity (recall $\v_{i,k+1}\in\myset{X}_{k-1}$)
\begin{multline*}
-(\tilde{\g}_{i,k}-\g_{i,k}^*)^\transp (\v_{i,k+1}-  \y_{i,k}^{*}) = \\ 
\tilde{\g}_{i,k}^\transp (\y_{i,k}^{*}-\v_{i,k+1}) + \g_{i,k}^{*,\transp}(\v_{i,k+1}-  \y_{i,k}^{*}) \leq \\  \myexp_{\myom,k}[f_{i,k}(\y_{i,k}^{*},\myom)- f_{i,k}(\v_{i,k+1},\myom)] + \g_{i,k}^{*,\transp}(\v_{i,k+1}-\y_{i,k}^{*}) + \epsilon \\ \leq -\frac{m_f}{2} \|\v_{i,k+1}-  \y_{i,k}^{*}\|^2 + \epsilon.
\end{multline*}
Combining these results with~\eqref{lemma.eq2} and \eqref{lemma.eq1bis} we obtain,  
\begin{multline*}
\|\y_{i,k+1} - \y_{i,k}^{*}\|^2 \leq (1+\mu)\Big[\|\v_{i,k+1}-  \y_{i,k}^{*}\|^2 +\\
 \alpha^2\Big(\frac{\epsilon^2}{4 |X_{k,k-1}|^2} + 2 L \epsilon + L^2 \|\v_{i,k+1}-\y_{i,k}^*\|^2\Big) +\\ \alpha\Big(-m_f \|\v_{i,k+1}-  \y_{i,k}^{*}\|^2 + 2 \epsilon\Big)\Big] + (1+1/\mu)\alpha^2 G^2, 
\end{multline*}
and therefore the claim is proven.
\end{myproof}
%------------------------------------------------------------------------------------
%------------------------------------------------------------------------------------
%------------------------------------------------------------------------------------

%and
%\begin{equation*}
%\|\y_{k+1} - \y_{k}^{*}\|^2 \leq \varphi(\myalp_k) \left\|\v_{k+1}-  \y_{k}^{*}\right\|^2 + \psi(\myalp_k, \myeps_k).
%\end{equation*}

We turn now our attention to the term $\|\v_{i,k+1} - \y_{i,k}^*\|^2$, which can be bounded as follows. 

\vskip0.1cm
\noindent\colorbox{mygray}{\begin{minipage}{0.48\textwidth}
\begin{lemma}\label{lemma.3}
Let $\y_{k}$ be the stacked version of the local decision variables $\y_{i,k}$, and let $\bar{\W} = \myexp[\W_k]$. Under the same assumptions and with the same definitions of Lemma~\ref{lemma.2} we have
\begin{equation*}
\myexp\left[\|\y_{k+1} - \y_{k}^{*}\|^2\right] \hskip-0.05cm\leq \hskip-0.05cm\varphi(\alpha,\bar{\W},\mu) \myexp\left[\|\y_{k}-  \y_{k}^{*}\|^2\right] + \alpha \psi(\alpha, \epsilon,\mu),
\end{equation*}
\end{lemma}\end{minipage}}
\noindent\colorbox{mygray}{\begin{minipage}{0.48\textwidth}
where 
\begin{eqnarray*}
\varphi(\alpha, \bar{\W},\mu) &=& (1+\mu)(1+\alpha^2L^2- \alpha m_f)(1 - \beta \lambda_2(\bar{\W}))\\
\psi(\alpha, \epsilon,\mu) &=&  (1+\mu)n\epsilon\Big(\alpha\frac{\epsilon}{4 |\myset{X}_{k,k-1}|^2} + 2 \alpha L  + 2 \Big) + \\ &&\qquad\qquad\qquad (1+1/\mu)\alpha n G^2.
\end{eqnarray*}
%\end{lemma}
\end{minipage}}

\vskip0.1cm
\noindent\begin{myproof}
By the definition of $\v_{i,k+1}$ in~\eqref{v}, we can write its stacked version for all $i$'s as
\begin{equation*}
\v_{k+1} = \y_{k} - \beta  (\W_{k}\otimes {\bf I}_d) \y_{k},
\end{equation*}
and since $ (\W_{k}\otimes {\bf I}_d) \y_{k}^* = {\bf 0}_{nd}$, we can expand 
\begin{equation*}
\|\v_{k+1} - \y_{k}^*\|^2 = (\y_k - \y_k^*)^\transp({\bf I}_{nd} - \beta \W_k \otimes {\bf I}_d)^2 (\y_k - \y_k^*).
\end{equation*}
We use now the fact that $\beta < 1/n$, which also means\footnote{In fact, $\lambda_1({\bf I}_{nd} - \beta {\W}_k\otimes {\bf I}_d)= 1 - \beta \lambda_n(\W_k)$; bounding $\lambda_n(\W_k)\leq n $ from~\cite{Rojo2007}, the result follows.} ${\bf I}_{nd} \succeq {\bf I}_{nd} - \beta {\W}_k\otimes {\bf I}_d \succ 0$, therefore  ${\bf I}_{nd} - \beta {\W}_k\otimes {\bf I}_d \succeq ({\bf I}_{nd} - \beta {\W}_k\otimes {\bf I}_d)^2$, and thus
\begin{multline*}
\|\v_{k+1} - \y_{k}^*\|^2 = (\y_k - \y_k^*)^\transp({\bf I}_{nd} - \beta \W_k \otimes {\bf I}_d)^2 (\y_k - \y_k^*) \leq \\
(\y_k - \y_k^*)^\transp({\bf I}_{nd} - \beta \W_k \otimes {\bf I}_d) (\y_k - \y_k^*).
\end{multline*}
We take now the expectation of the previous expression  
\begin{multline*}
\mathbb{E}\left[ \|\v_{k+1} - \y_{k}^*\|^2\right] \leq \\ \mathbb{E}\left[(\y_k - \y_k^*)^\transp ({\bf I}_{nd} - \beta \W_k \otimes {\bf I}_d)(\y_k - \y_k^*)\right],
\end{multline*}
and since $\y_k$ and $\W_k$ are independent, 
\begin{multline*}
\mathbb{E}\left[ \|\v_{k+1} - \y_{k}^*\|^2\right] \leq\\ \mathrm{tr}\left( ({\bf I}_{nd} - \beta \bar{\W} \otimes {\bf I}_d)\mathbb{E}\left[ (\y_{k} - \y_{k}^*)(\y_{k} - \y_{k}^*)^\transp\right]\right),
\end{multline*}
where $\mathrm{tr}(\cdot)$ is the trace operator. 
 And finally, given that $\bar{\W} (\y_k - \y_k^*)$ will be zero only when we reach consensus (i.e., at optimality), 
\begin{multline*}
\mathrm{tr}\left( ({\bf I}_{nd} - \beta \bar{\W} \otimes {\bf I}_d)\mathbb{E}\left[ (\y_{k} - \y_{k}^*)(\y_{k} - \y_{k}^*)^\transp\right]\right)\leq \\ (1 - \beta \lambda_2(\bar{\W}))\myexp\left[\|\y_k - \y_k^*\|^2\right],
\end{multline*}
and the claim follows by combining this with Lemma~\ref{lemma.2}:
\begin{multline*}
\mathbb{E}\Big[\sum_{i\in V} \|\y_{i,k+1} - \y_{i,k}^{*}\|^2\Big] \leq \\ \mathbb{E}\Big[\sum_{i\in V} (1+\mu)(1+\alpha^2L^2- \alpha m_f)\left\|\v_{i,k+1}-  \y_{i,k}^{*}\right\|^2\Big] +\\ \alpha\psi(\alpha,\epsilon,\mu)  \leq \varphi(\alpha, \bar{\W},\mu) \mathbb{E}\Big[\|\y_{k} - \y_{k}^*\|^2\Big] + \alpha\psi(\alpha,\epsilon,\mu).
\end{multline*}
\end{myproof}

\subsection{Proof of Theorem~\ref{convergence}}

We are now ready to prove Theorem~\ref{convergence}. 

\noindent\begin{myproof}
We use again the Peter-Paul inequality, Assumption~\ref{as.opt}, and Lemma~\ref{lemma.3} to establish
\begin{multline*}
\mathbb{E}\left[\|\y_{k+1} - \y_{k+1}^{*}\|^2\right] \leq (1+\mu)\mathbb{E}\left[ \|\y_{k+1} - \y_{k}^{*}\|^2\right] + \\(1+1/\mu)\mathbb{E}\left[\|\y_{k+1}^* - \y_{k}^{*}\|^2\right] \leq n \delta_\x^2 (1+1/\mu) + \\ \left(\varphi(\alpha, \bar{\W},\mu)  \mathbb{E}\left[\|\y_{k}-  \y_{k}^{*}\|^2 \right]+ \alpha\psi(\alpha, \epsilon,\mu)\right)(1+\mu),%\\ + n \delta_\x^2 (1+1/\mu), 
\end{multline*}
which is true for any scalar $\mu>0$. In order to ensure convergence of the sequence $\{\y_k\}$, it has to be 
\begin{equation*}
r := \varphi(\alpha, \bar{\W},\mu)(1+\mu) < 1, %\quad \textrm{so that}
\end{equation*}
so that, 
\begin{multline*}
\mathbb{E}\left[\|\y_{k+1} - \y_{k+1}^{*}\|^2\right] \leq  r^{k}  \mathbb{E}\left[\|\y_{1} - \y_{1}^{*}\|^2\right]  + \\ \sum_{\tau = 1}^{k}  r^{\tau}\left(\alpha\psi(\alpha, \epsilon,\mu)(1+\mu) + n \delta_\x^2 (1+1/\mu)\right), 
\end{multline*}
and finally, $\liminf_{k\to \infty}\mathbb{E}\left[\|\y_{k+1} - \y_{k+1}^{*}\|^2\right]  =$
\begin{multline*}
\liminf_{k\to \infty}\sum_{\tau = 1}^{k}  r^{\tau}\left(\alpha\psi(\alpha, \epsilon)(1+\mu) + n \delta_\x^2 (1+1/\mu)\right)  \\ \leq \frac{1}{1-r} \Big(\alpha\psi(\alpha, \epsilon,\mu)(1+\mu) + n \delta_\x^2 (1+1/\mu)\Big).
\end{multline*}
To ensure $r<1$, we need to have
\begin{equation*}%\label{alpha}
(1 + \alpha^2 L^2- \alpha m_f) (1 - \beta \lambda_2(\bar{\W})) (1+\mu)^2 < 1.
\end{equation*}
We call $\gamma = 1 - \beta \lambda_2(\bar{\W})$, and by the definition of $\beta$ and Assumption~\ref{as.persistent}, it is $0<\gamma<1$.
%Choose now $\alpha < m_f/L^2$ and $\mu =  (m_f/L^2)/(\alpha \gamma) - 1$
%
%
%
%
%The condition~\eqref{alpha} implies
%\begin{equation*}\label{alpha1}
%\alpha \in \left(\frac{m_f - \sqrt{\triangle}}{4L^2} , \frac{m_f + \sqrt{\triangle}}{4L^2}\right),
%\end{equation*}
%with
%$$
%\triangle = m_f^2 + 4L^2 \left(\frac{1}{\gamma(\mu+1)} -1 \right).
%$$
%For $\alpha$ to have possible values, $\triangle$ has to be positive; given the fact that $m_f\leq L$ (see Assumption~\ref{as.function}), then a sufficient condition for $\triangle > 0$ is $\mu < \frac{4}{3\gamma} -1$. A loser sufficient condition (but independent of $\gamma$) is therefore $\mu< \nicefrac{1}{3}$. This proves the first part of the theorem and the linear convergence rate. 
%
%To prove the second part we notice that for $\triangle$ we have two distinct cases, either $\triangle \geq m_f^2$, and in this case we have only an upper bound on $\alpha$, or  $0< \triangle < m_f^2$, which yields the existence of a lower bound. The case $\triangle \geq m_f^2$ is attained for $\gamma \leq 1/(\mu +1)$, meaning $\beta \lambda_2(\bar{\W}) \geq \mu/(\mu+1)$, and therefore the claim. 
By choosing $\mu = 1/\sqrt{\gamma}-1$, we obtain the condition
%\begin{equation*}
$ \varrho := 1 + \alpha^2 L^2- \alpha m_f  < 1,$ 
%\end{equation*}
that is $\alpha < m_f/L^2$, and the limit result
\begin{multline*}
\liminf_{k\to \infty}\myexp[\|\y_{k+1} - \y_{k+1}^{*}\|^2]\leq  \\
\frac{1}{1-\varrho} \Big(\psi(\alpha, \epsilon,\mu)\frac{\alpha}{\sqrt{\gamma}} + n\delta_\x^2 \frac{1}{1-\sqrt{\gamma}}\Big),
\end{multline*}
with $\mu = 1/\sqrt{\gamma} - 1$, and the theorem is proven. 
\end{myproof}

\subsection{Proof of Theorem~\ref{policy}}

%The role of the policy $\pi$ in Theorem~\ref{policy} is to keep the error term $\psi(\alpha, \myeps_\tau)$ bounded by a specified threshold (and thereby guaranteeing that $\tilde{\g}_{i,k}$ is indeed an $\epsilon$-gradient). If we indicate with $\myom'$ the available (and possibly outdated) information at node $i$, then

\begin{myproof}
We start expanding the condition~\eqref{cond} by using the triangle inequality as 
\begin{multline*}   
\hskip-.4cm \left\|\myexp_{\myom_j\hskip-0.05cm,k} [\nabla_{\x} f_{i,k,j}(\x,\myom_j)|_{\v_{i,k+1}}\hskip-0.05cm] - \myexp_{\myom_j'}\hskip-0.05cm[\nabla_{\x} f_{i,k,j}(\x,\myom_j)|_{\v_{i,k+1}}\hskip-0.05cm]\right\| \\ \leq
\left\|\myexp_{\myom_j,k} [\nabla_{\x} f_{i,j}'(\x,\myom_j)|_{\v_{i}'}] - \myexp_{\myom_j'} [\nabla_{\x} f_{i,j}'(\x,\myom_j)|_{\v_{i}'}] \right\| + \\ \left\|\myexp_{\myom_j,k} \big[\nabla_{\x} f_{i,k,j}(\x,\myom_j)|_{\v_{i,k+1}}-\nabla_{\x} f_{i,j}'(\x,\myom_j)|_{\v_{i}'}\big]\right.\\-\left.\myexp_{\myom_j'} \big[\nabla_{\x} f_{i,k,j}(\x,\myom_j)|_{\v_{i,k+1}}-\nabla_{\x} f_{i,j}'(\x,\myom_j)|_{\v_{i}'}\big]\right\|.
\end{multline*}
The right-most term can be bounded as
\begin{multline*}   
\left\|\myexp_{\myom_j,k} \big[\nabla_{\x} f_{i,k,j}(\x,\myom_j)|_{\v_{i,k+1}}-\nabla_{\x} f_{i,j}'(\x,\myom_j)|_{\v_{i}'}\big]\right.\\-\left.\myexp_{\myom_j'} \big[\nabla_{\x} f_{i,k,j}(\x,\myom_j)|_{\v_{i,k+1}}-\nabla_{\x} f_{i,j}'(\x,\myom_j)|_{\v_{i}'}\big]\right\| \leq \\ \int_{\Omega_j} \Big\|\nabla_{\x} f_{i,k,j}(\x,\myom_j)|_{\v_{i,k+1}}-\nabla_{\x} f_{i,j}'(\x,\myom_j)|_{\v_{i}'}\Big\| \, \cdot \\ \Big|p_{\myom_{j},k}(\myom)-p_{\myom_{j}'}(\myom)\Big| \mathrm{d}\myom \leq   \nu \mathcal{U}_{\mathrm{R},ij}(\nabla_{\x} f_{i,k,j}|\nabla_{\x} f_{i,j}'),
\end{multline*}
where we have used the policy $\pi_{ij}$ to bound $|p_{\myom_{j},k}(\myom)-p_{\myom_{j}'}(\myom)|$. Therefore we can write
\begin{multline*}   
\hskip-.4cm \left\|\myexp_{\myom_j\hskip-0.05cm,k} [\nabla_{\x} f_{i,k,j}(\x,\myom_j)|_{\v_{i,k+1}}\hskip-0.05cm] - \myexp_{\myom_j'}\hskip-0.05cm[\nabla_{\x} f_{i,k,j}(\x,\myom_j)|_{\v_{i,k+1}}\hskip-0.05cm]\right\| \\ \leq \mathcal{U}_{\mathrm{S},ji}^{(1)}(\myom_k|\myom') + \nu \mathcal{U}_{\mathrm{R},ij}(\nabla_{\x} f_{i,k,j}|\nabla_{\x} f_{i,j}').
\end{multline*}
And finally, a sufficient condition for~\eqref{cond} is 
\begin{eqnarray*}
&\mathcal{U}_{\mathrm{S},ji}^{(1)}(\myom_k|\myom') \leq  (1-\eta) \frac{\epsilon}{2 |X_{k,k-1}| \Delta_i}, \\  &\mathcal{U}_{\mathrm{R},ij}(\nabla_{\x} f_{i,k,j}|\nabla_{\x} f_{i,j}') \leq \frac{\eta}{\nu}\frac{\epsilon}{2 |X_{k,k-1}| \Delta_i}
\end{eqnarray*}
for any $\eta\in[0,1]$; this together with the bound on $\mathcal{U}_{\mathrm{S},j}^{(2)}(\myom_k|\myom')$ yields the claim. 
\end{myproof}

\section{Conclusions}

We have proposed a distributed stochastic $\epsilon$-gradient asynchronous algorithm to optimize a rather general convex separable time-varying stochastic program. To alleviate the possibly high communication demands among the nodes to track variations in the \textsc{pdf}'s of the stochastic variables, we have devised a utility-based policy to trigger the said communication. The overall scheme converges linearly to an error bound whose size depends on the constant stepsize $\alpha$, on $\epsilon$, and the variability in time of the optimizer. Initial simulation results are encouraging and well display the added value of the proposed approach.   

Many relevant questions require further studies. First of all, the Peter-Paul inequality is known to be a loose bound; by substituting this inequality with tighter bounds, we expect more accurate asymptotic results. The algorithm should be further modified to allow for less restrictive assumptions, both on the cost functions and on the stochastic variables. Finally, in order to be able to characterize dynamical constraints, it will be interesting to explicitly include in the analysis cost functions that depend on past variables generated locally, e.g., $f_{i,k}$ should also depend on $\y_{i,k-
1}$, as pointed out in Remark~\ref{remdin1}.

 \footnotesize
\bibliographystyle{ieeetran}
\bibliography{../../PaperCollection2}
\end{document}